\documentclass{siamltex}
\usepackage{amsfonts,amsmath,mathrsfs}
\usepackage{amssymb}
\usepackage{graphicx}
\usepackage{wrapfig}
\usepackage{subfig}
\usepackage{float}
\usepackage{bbm}

\usepackage{layouts}

\newcommand{\N}{\mathbb{N}}
\newcommand{\R}{\mathbb{R}}

\newtheorem{remark}[theorem]{Remark}

\DeclareMathOperator{\isoTD}{\mathtt{isoTD}}

\newcommand{\domain}{D}
\newcommand{\stokastinendomain}{\Omega}

\newcommand{\dd}{\textrm{d}}
\newcommand{\y}{\textbf{y}}
\newcommand{\x}{\textbf{x}}
\newcommand{\vbf}{\mathbf{v}}

\newcommand{\vv}{\mathrm{v}}
\newcommand{\transposeT}{\mathsf{T}}

\title{Stochastic Galerkin finite element method with local conductivity basis for electrical impedance tomography}

\author{
N.~Hyv\"onen\footnotemark[2]
\and M.~Leinonen\footnotemark[2]
}

\begin{document}
\maketitle

\renewcommand{\thefootnote}{\fnsymbol{footnote}}

\footnotetext[2]{
Aalto University, Department of Mathematics and Systems Analysis, P.O. Box 11100, FI-00076 Aalto, Finland (nuutti.hyvonen@aalto.fi, matti.leinonen@aalto.fi). This work was supported by the Academy of Finland (decision 267789) and the Finnish Doctoral Programme in Computational Sciences FICS.
}

\begin{abstract}
The objective of electrical impedance tomography is to deduce information about the conductivity inside a physical body from electrode measurements of current and voltage at the object boundary. In this work, the unknown conductivity is modeled as a random field parametrized by its values at a set of pixels. The uncertainty in the pixel values is propagated to the electrode measurements by numerically solving the forward problem of impedance tomography by a stochastic Galerkin finite element method in the framework of the complete electrode model. For a given set of electrode measurements, the stochastic forward solution is employed in approximately parametrizing the posterior probability density of the conductivity and contact resistances. Subsequently, the conductivity is reconstructed by computing the {\em maximum a posteriori} and {\em conditional mean} estimates as well as the posterior covariance. The functionality of this approach is demonstrated with experimental water tank data.
\end{abstract}

\renewcommand{\thefootnote}{\arabic{footnote}}



\pagestyle{myheadings}
\thispagestyle{plain}
\markboth{N.~HYV\"ONEN AND M.~LEINONEN}{SGFEM WITH LOCAL CONDUCTIVITY BASIS FOR EIT}

\section{Introduction}
The aim of {\em electrical impedance tomography} (EIT) is to retrieve useful information about the conductivity inside an examined physical body based on boundary measurements of current and voltage. In practice, the boundary data are gathered with a finite number of contact electrodes; the most accurate model for EIT is the {\em complete electrode model} (CEM) \cite{Cheng89,Somersalo92}, which takes into account the electrode shapes and the contact resistances at the electrode-object interfaces. EIT has potential applications in, e.g., medical imaging, monitoring of industrial processes, and nondestructive testing of materials; see the review articles \cite{Adler11,Borcea02,Cheney99,Lionheart03,Uhlmann09} and the references therein for more information on EIT and related mathematics.

This work considers EIT from the standpoint of uncertainty quantification. The to-be-reconstructed conductivity is modeled as a random field parametrized by uniformly distributed mutually independent random variables representing the conductivity levels at a set of pixels. The range of the pixel values is chosen based on prior information, while the number of pixels is mainly dictated by computational restrictions. The contact conductances, i.e., the reciprocals of the contact resistances, are also assigned uniform prior densities. For a given measurement configuration, the uncertainty in the conductivity field and the contact resistances is propagated to the electrode measurements by approximately solving the stochastic version of the CEM forward problem by a {\em stochastic Galerkin finite element method} (sGFEM) \cite{Ghanem03,Schwab11a}, which in our case corresponds to discretizing the spatial domain by piecewise linear FEM basis functions and the stochastic domain by a spectral Galerkin method with a Legendre polynomial basis (cf.~\cite{Xiu10}). These steps can be carried out off-line,~i.e.,~prior to the actual measurements, assuming the measurement geometry as well as the ranges for the conductivity and contact conductance values are known in advance. 

After the electrode potentials corresponding to a set of applied current patterns have been measured, the stochastic forward solution can be used to explicitly write an approximate parametrization for the {\em posterior} density of the conductivity,~i.e.,~for the {\em posterior} of the pixelwise conductivity levels. At this stage, it is also possible to `update' the prior in case one has more specific information on the particular conductivity at hand. 
In this work, the information on the range of the pixelwise conductivity levels assumed in the forward solver is complemented by a Gaussian smoothness prior, but we want to emphasize that other forms of {\em a priori} information could as well be incorporated in the inverse solver. The actual conductivity reconstructions are obtained by computing {\em maximum a posteriori} (MAP) and {\em conditional mean} (CM) estimates,~i.e.,~the maximum point and the expected value of the approximate posterior density, respectively. In our setting, the computation of the former corresponds to minimizing a high-dimensional positive-valued polynomial, whereas the latter deals with high-dimensional integration with an explicitly known integrand. The reconstructions of the conductivity are complemented with visualizations of the posterior standard deviation. 

The papers \cite{Leinonen14,Hakula14} introduced a reconstruction method for two-dimensional EIT by applying sGFEM to the CEM under the assumption that the conductivity is {\em a priori} known to be a lognormal random field. To be more precise, the conductivity was parametrized using its truncated exponential Karhunen--Lo\`eve expansion, and reconstructions were computed by estimating the random coefficients in the truncated expansion on the basis of (simulated) measurement data. Although the assumption of lognormality can be considered natural \cite{Leinonen14},
the major drawback of the approach in \cite{Hakula14} is that
the spatial and stochastic components of the sGFEM solution cannot be decoupled, which results in relatively full system matrices (cf.~\cite[Section~6.1]{Hakula14}). This can easily be a deal-breaker in practical EIT since the accurate enough solution of the stochastic CEM forward model by sGFEM requires the use of a {\em high} number of degrees of freedom. 
The algorithm presented in this work can be considered a modified version of the one in \cite{Hakula14}, aiming at better computational feasibility: The pixelwise parametrization by uniformly distributed random variables results in a very sparse sGFEM system and it also allows trivial control over the positivity of the conductivity. Compared with \cite{Hakula14}, our new algorithm makes it possible to straightforwardly update the prior information on the conductivity in the on-line solution phase and to consider the estimation of a higher number of parameters from electrode measurements, resulting in improved reconstructions.

Compared with previous Bayesian techniques for tackling the inverse problem of practical EIT (see,~e.g.,~\cite{Darde13b,Heikkinen02,Kaipio00,Karhunen10} and the references therein), the main advantage of our approach is the following: Our method produces an (approximate) parametrization of the posterior density, i.e.,~of the idealized solution to the inverse problem in the Bayesian sense, which makes it possible to analyze the posterior without referring to the elliptic boundary value problem associated to the CEM. (In the `standard' Bayesian approach to EIT, each evaluation of the posterior density requires solving as many deterministic CEM forward problems as there are applied current patterns.) In particular, if the sGFEM solution of the CEM has been computed prior to the measurements, reconstructions and corresponding uncertainty estimates for the conductivity can be produced without ever returning to the CEM forward problem itself. This leads to obvious computational benefits because evaluating explicitly known functions is typically cheaper than solving several elliptic boundary value problems. The obvious disadvantage of the proposed method is the requirement of precomputing an accurate enough  sGFEM forward solution for the CEM. However, the inevitable increase in computational resources and further development of stochastic finite element algorithms (see,~e.g.,~\cite{Bieri09a}) may well facilitate a satisfactory solution to this problem in the future.

The approach of this work is purely computational: based on experimental data from water tank experiments, we demonstrate that the introduced algorithm produces two-dimensional reconstructions that are arguably almost as good as the state-of-the-art Bayesian reconstructions from experimental data under a {\em smoothness prior} (cf.,~e.g.,~\cite{Darde13b,Karhunen10}). For information on the convergence of the sGFEM-parametrized posterior density in closely related settings, we refer to \cite{Schillings14,Schwab12} and the references therein. However, we are not aware of proper convergence analysis of sGFEM-based reconstruction algorithms for inverse elliptic {\em boundary} value problems. Moreover, to the best of our knowledge, this is the first time that any stochastic finite element method has been employed to compute EIT reconstructions from experimental data. See \cite{Dashti13,Dashti11,Hoang13,Schillings13,Schillings14,Simon14,Stuart10} for related approaches to solving inverse problems.

The rest of this paper is organized as follows.
The {\em stochastic complete electrode model} (SCEM) is introduced in Section~\ref{sec:SCEM}, and solving the SCEM forward problem by sGFEM is considered in Section~\ref{sec:forward_solution}.  We focus on the Bayesian inverse problem of EIT in Section~\ref{sec:inverse_solution}, and Section~\ref{sec:implementation} discusses the two-phase implementation of our reconstruction algorithm. The numerical examples are presented in Section~\ref{sec:numerical_examples}. We conclude with a few remarks in Section~\ref{sec:conclusions}.

\section{Stochastic complete electrode model}
\label{sec:SCEM}
In this section, we introduce the SCEM for modeling practical EIT measurements with a random conductivity and contact resistances. For the traditional deterministic formulation together with its physical and experimental justification, see \cite{Cheng89,Somersalo92}.

Let $\domain\subset\R^n$, $n=2$ or $3$, be a bounded domain with a smooth enough boundary and let $(\stokastinendomain,\Sigma,P)$ be a probability space. We interpret the internal conductivity of $D$ as a random field $\sigma(\cdot,\cdot): \stokastinendomain \times \domain \rightarrow \R$ which is assumed to be a uniformly strictly positive element of $L^\infty(\stokastinendomain \times D)$, i.e.,
\begin{align*}
P\left(\omega \in \stokastinendomain \ : \ \sigma_{\textrm{min}} \leq 
\operatorname*{ess\,inf}_{\x \in \domain}
\sigma(\omega,\x) \leq
\operatorname*{ess\,sup}_{\x \in \domain}
\sigma(\omega,\x) \leq  \sigma_{\textrm{max}}\right) = 1
\end{align*}
for some constants $\sigma_{\textrm{min}}, \sigma_{\textrm{max}} > 0$.
The perfectly conducting electrodes $E_1, \dots, E_M$, $M \in \N \setminus \{1\}$, attached to $\domain$ are identified with the corresponding open, connected, and mutually disjoint subsets of $\partial \domain$. We denote $E = \cup_m E_m$, $I = [I_1, \dots , I_M]^{\transposeT}$, and $U = [U_1, \dots, U_M]^{\transposeT}$, where $I_m\in\R$ and $U_m: \Omega \to \R$ are the injected deterministic net current and the measured random voltage, respectively, on the $m$th electrode. 
The current pattern $I$ belongs to the mean-free subspace $\R^{M}_\diamond$ of $\R^{M}$ by virtue of the conservation of charge; the voltage vector $U$ is interpreted as a (random) element of  $\R^{M}_\diamond$ by choosing the ground level of potential appropriately.
The contact resistances representing the resistive layers between the electrodes and the domain $D$ are modeled by random variables $z_m:\stokastinendomain \rightarrow \R$, $m=1,\ldots,M$, which are assumed to be uniformly strictly positive and bounded: 
\begin{align*}
P(\omega \in \stokastinendomain \ : \ z_{\textrm{min}} \leq z_m(\omega) \leq z_{\textrm{max}}) = 1,\qquad m=1,\ldots,M, 
\end{align*}
for some $z_{\textrm{min}}, z_{\textrm{max}} > 0$.

Denote $\mathcal{H} := H^1(D) \oplus \R^M_\diamond$ and let us introduce the Bochner space
\[
L_P^2(\stokastinendomain;\mathcal{H}) := \left\{ (u,U):\stokastinendomain\rightarrow \mathcal{H} \ \big| \  \int_\stokastinendomain \|(u(\omega),U(\omega))\|_{\mathcal{H}}^2 \, \dd P(\omega)<\infty \right\}
\]
that allows the decomposition $L_P^2(\stokastinendomain;\mathcal{H})\simeq L_P^2(\stokastinendomain)\otimes \mathcal{H}$, where $\otimes$ denotes the tensor product between Hilbert spaces (cf., e.g.,~\cite{Schwab11a}).
The SCEM forward problem is as follows. For a given deterministic electrode current pattern $I \in \R^{M}_\diamond$, find a pair $(u,U) \in L_P^2(\stokastinendomain;\mathcal{H})$ that satisfies the following boundary value problem $P$-almost surely:
\begin{align*}
\begin{array}{ll}
\nabla \cdot (\sigma \nabla u) = 0 \qquad  &\text{in} \ \domain, \\[8pt] 
{\displaystyle \frac{\partial u}{\partial \nu}} = 0 \qquad &\text{on} \ \partial \domain\setminus\overline{E},\\[2mm] 
{\displaystyle u+z_m \sigma \frac{\partial u}{\partial \nu}}= U_m \qquad &\text{on} \ E_m, \quad m=1, \dots, M, \\[3mm] 
{\displaystyle \int_{E_m} \sigma \frac{\partial u}{\partial \nu} \, \dd S} = I_m, \qquad & m=1,\ldots,M, 
\end{array}
\end{align*}
where $\nu = \nu(x)$ is the exterior unit normal of $\partial D$.
The corresponding variational formulation is to find $(u,U) \in L_P^2(\stokastinendomain;\mathcal{H})$ such that
\begin{align}
\label{equ:stokamuoto}
\mathbb{E} \big[ B\big((u,U),(v,V)\big) \big] \, = \, 
I \cdot \mathbb{E}[V] \qquad \textrm{for} \ \textrm{all} \ (v,V) \in 
L_P^2(\stokastinendomain;\mathcal{H}),
\end{align}
where $\mathbb{E}[\, \cdot \,]$ denotes the expectation
and the bilinear form 
$B:  \mathcal{H} \times \mathcal{H} \to \R$ is defined via 
\[
B\big((u,U),(v,V)\big) \, = \, \int_{\domain} \sigma \nabla u \cdot \nabla v \,\dd \x + 
\sum_{m=1}^{M} \frac{1}{z_m} \int_{E_m}( U_m -u ) (V_m - v )\, \dd S. 
\]
The unique solvability of the SCEM forward problem can be proved by extending the deterministic argumentation in~\cite{Somersalo92}.

\subsection{Parametric deterministic SCEM}
\label{sec:paramdetSCEM}
In the rest of this work, the conductivity is assumed to be parametrized by its random values at a finite set of open pixels $D_1, \dots, D_L$, which constitute a partition of $\domain$, i.e., $\overline{D} = \cup \overline{D}_l$. More precisely,
\begin{align}
\label{equ:klexpansion0}
\sigma(\omega,\x) = \sigma_0 + \sum_{l=1}^L \sigma_l \mathbf{1}_{D_l}(\x) Y_l(\omega), \qquad \omega \in \Omega, \ \x \in D,
\end{align}
 where 
$\sigma_0 \in \R_+$, $\sigma_l \in \R_+ \cup \{ 0 \}$, and $ \sigma_l < \sigma_0$ for $l=1,\ldots,L$. Moreover, 
$\mathbf{1}_{D_l}$ is the indicator function of $D_l$, and each random variable $Y_1, \dots , Y_L$ is uniformly distributed on the interval $[-1,1]$.
For every $m=1, \dots, M$, the contact resistance $z_m$ is assumed to follow the inverse uniform distribution on the interval $[b_m^{-1},a_m^{-1}]$, where $0 < a_m < b_m$.
In consequence, the contact conductances $\zeta_1:= z_1^{-1}, \dots, \zeta_M:=z_M^{-1}$ can be presented as 
\begin{align*}
\zeta_m(\omega) = \frac{1}{2}(a_m+b_m) + \frac{1}{2}(b_m-a_m)Y_{L+m}(\omega), \qquad m=1, \dots, M,
\end{align*}
where each $Y_{L+1}, \dots, Y_{L+M}$ obeys the uniform distribution on $[-1,1]$. 
It is assumed that $Y_1, \dots , Y_{L+M}$ are mutually independent. 

To simplify the notation, we define
$$
\mathbf{Y}_\sigma = (Y_1,\ldots,Y_L), \quad \mathbf{Y}_\zeta = (Y_{L+1},\ldots,Y_{L+M}),
$$
and denote $\mathbf{Y} = (\mathbf{Y}_\sigma, \mathbf{Y}_\zeta)$. In particular,
$\mathbf{Y}: \Omega \to \R^{L+M}$ has the probability density
\begin{align}
\label{equ:probden}
\rho(\mathbf{y}) =
\left\{
\begin{array}{ll}
2^{-(L+M)} &\qquad {\rm if} \ \mathbf{y} \in \Gamma, \\[1mm]
0 &\qquad {\rm otherwise},
\end{array}
\right.
\end{align}
where $\Gamma = [-1,1]^{L+M}$.

Substituting the above choices in \eqref{equ:stokamuoto}, we arrive at our parametric deterministic variational formulation of the SCEM forward problem: 
find $(u,U) \in L^2(\Gamma;\mathcal{H})$ such that
\begin{align}
\label{equ:paramuoto}
\int_{\Gamma} \! \Big[ \int_{\domain} \sigma(\y,\x) \nabla u \cdot \nabla v \,\dd \x  + \! \sum_{m=1}^{M} \zeta_m(\y) \int_{E_m}( U_m - u ) (V_m - v )  \, \dd S \Big] \dd \y
= 
I \cdot \! \int_\Gamma \! V(\y) \dd \y
\end{align}
for all $(v,V) \in L^2(\Gamma;\mathcal{H})$. Here, with a slight abuse of the notation,
\begin{align}
\label{equ:klexpansion}
\sigma(\y,\x) = \sigma_0 + \sum_{l=1}^L \sigma_l \mathbf{1}_{D_l}(\x) \, y_l
\end{align}
and
\begin{align}
\label{sconduct}
\zeta_m(\y) = \frac{1}{2}(a_m+b_m) + \frac{1}{2}(b_m-a_m)y_{L+m}, \qquad m=1, \dots, M,
\end{align}
i.e., we have interpreted the conductivity and the contact conductances as functions of the parameter vector $\y = (\y_\sigma, \y_\zeta) \in \Gamma \subset \R^{L+M}$. 

\begin{remark}
As the probability density \eqref{equ:probden} is piecewise constant,
we have dropped 
the `weight' $\rho(\mathbf{y})$ from the integrals in \eqref{equ:paramuoto} and refrained from introducing weighted $L^2$-spaces. In general, this is not recommendable; see, e.g., \cite{Schwab11a,Leinonen14}.
\end{remark}

\section{Stochastic forward solution}
\label{sec:forward_solution}
To numerically solve \eqref{equ:paramuoto}, we need to discretize
$L^2(\Gamma; \mathcal{H}) \simeq L^2(\Gamma) \otimes (H^1(\domain)\oplus \R^M_\diamond)$, which boils down to choosing finite-dimensional bases for (certain subspaces of) $L^2(\Gamma)$, $H^1(\domain)$, and $\R^M_\diamond$. The spaces $H^1(\domain)$ and $\R^M_\diamond$ are handled as in standard FEM, whereas for $L^2(\Gamma)$ we use the spectral Galerkin method with a multivariate Legendre polynomial basis. The latter choice is reasonable as \eqref{equ:paramuoto} includes no differentiation with respect to $\y$.
 
For $H^1(\domain)$ we use the standard FEM with piecewise linear basis $\{\varphi_j \}_{j=1}^{N_D} \subset H^1(\domain)$, $N_D \in \N$, with respect to a suitable mesh. As the mean-free basis vectors for $\R^M_\diamond$, we employ
\begin{equation}
\label{meanfreebasis}
\vv_{i}= {\mathrm e}_1 - {\mathrm e}_{i+1}, \quad i = 1, \dots, M-1,
\end{equation}
with ${\mathrm e}_i$ denoting the $i$th Euclidean basis vector of $\R^M$.
To introduce the discretization of $L^2(\Gamma)$, we first recall the definitions of the univariate and multivariate Legendre polynomials.

\begin{definition}[Legendre polynomials]
Let $m\in\N_0:= \N\cup \{0 \} = \{0,1,2,\ldots\}$. The $m$th {\em univariate Legendre polynomial} is defined as
\[
L_m(y) := \frac{\sqrt{2m+1}}{2^{m+1/2}\,m!}\frac{\dd^m}{\dd y^m}[(y^2-1)^m],
\]
where $y \in \R$.
\end{definition}

Note that we have (nonstandardly) normalized the Legendre polynomials so that they are orthonormal with respect to the $L^2$ inner product over $[-1,1]$:
\[
\int_{-1}^1 L_k(y) L_l(y) \, \dd y  =  \delta_{k,l}, \qquad k,l \in \N_0,
\]
where $\delta_{k,l}$ is the Kronecker's delta.

\begin{definition}[Multivariate Legendre polynomials]
\label{def:multivariatepolynomial}
Let $P \in \N$ and $\mu \in \N_{0}^P$
be a multi-index. The {\em multivariate Legendre polynomial} $L_\mu$, also called {\em chaos polynomial}, is defined as 
\begin{align*}
L_\mu(\mathbf{y}) := 
\prod_{k=1}^{P} L_{\mu_k}(y_k), \qquad \mathbf{y} \in \R^P,
\end{align*}
where $L_{\mu_k}$ is the $\mu_k$th univariate Legendre polynomial. 
\end{definition}

The set $\mathcal{P} :=\{ L_{\mu}~|~\mu\in \N_{0}^{L+M}\}$
is an orthonormal basis of $L^2(\Gamma)$ (cf.,~e.g.,~\cite{Schwab11a}),
and thus any function $f \in L^2(\Gamma)$ admits a {\em polynomial chaos} representation, 
\begin{align}
\label{equ:chaos_rep}
f \, =  \!\! \sum_{\mu \in \N_{0}^{L+M}}\!\! \big( f,L_\mu \big)_{L^2(\Gamma)} L_\mu 
\end{align}
in the topology of $L^2(\Gamma)$.
In practical computations the number of multi-indices considered in \eqref{equ:chaos_rep} must naturally be finite, and hence we must replace $\N_{0}^{L+M}$ with a finite subset of multi-indices $\Lambda \subset \N_{0}^{L+M}$.

The set $\Lambda$ is ideally chosen so that
\begin{align*}
f  \, \approx \, \sum_{\mu \in \Lambda}\!\! \big( f,L_\mu \big)_{L^2(\Gamma)} L_\mu
\end{align*}
is as accurate as possible for the considered $f$ under a given constraint on the cardinality $\# \Lambda$. When solving \eqref{equ:paramuoto}, one would like to get good representations (for the FEM approximations) of $f = u(\, \cdot \, , \, \x )$, $\x \in D$. In practice, estimating {\em a priori} optimal index sets for the solutions of \eqref{equ:paramuoto} is highly nontrivial (but possible to a certain extent~\cite{Bieri09a}), and hence we resort in this work to generic index sets which are easy to generate and give equal weight to each dimension in $\Gamma$.

\begin{definition}[Isotropic total degree index set]
Let $P, Q \in \N$. The 
$\isoTD$ index set is defined as
\begin{align*}
\isoTD(P,Q) = \left\{\mu \in  \N_{0}^{P}~\big|~\sum_{k=1}^{P} \mu_k \leq Q \right\}.
\end{align*}
\end{definition}
It is easy to see that the cardinality of the $\isoTD(P,Q)$ index set is
\begin{equation}
\label{stoch_df}
\# \isoTD(P,Q) \, = \,  {P+Q \choose Q} \, .
\end{equation}
In what follows, we use $\Lambda = \isoTD(L+M,Q)$ for some $Q \in \N$ and denote $N_\Gamma = \# \Lambda$. See, e.g., \cite{Back11,Beck12,Bieri09a} and the references therein for information on other types of index sets.

We look for an approximation $(\tilde{u},\tilde{U})$ of the parametric deterministic SCEM solution $(u,U)$ to \eqref{equ:paramuoto} in the form
\begin{subequations}
\begin{align}
\label{sGFEM_u}
u(\y,\x)\approx\tilde{u}(\y,\x) &=
\sum_{j=1}^{N_D}
\sum_{\mu \in\Lambda}
\alpha_{j,\mu}L_{\mu}(\y)\varphi_j(\x),\\
\label{sGFEM_U}
U(\y)\approx\tilde{U}(\y) &= \sum_{i=1}^{M-1} \sum_{\mu\in\Lambda}\beta_{i,\mu} L_{\mu}(\y) \vv_{i},
\end{align}
\end{subequations}
where $\{\alpha_{j,\mu}\}\subset\R$ and $\{\beta_{i,\mu}\}\subset\R$ are the to-be-determined real coefficients.
In particular, the approximation of the electrode potentials in \eqref{sGFEM_U} is an $M$-dimensional vector whose components are $Q$th order polynomials in $\y$. We denote by  
$\alpha \in \R^{N_D N_\Gamma}$ and $\beta \in \R^{(M-1)  N_\Gamma}$ the block vectors
defined by $\{\alpha_{j}\}_{\mu} = \alpha_{j,\mu}$ 
and  $\{\beta_{i}\}_{\mu} = \beta_{i,\mu}$, respectively.

The coefficient vector $(\alpha, \beta)$ is determined via the standard Galerkin projection:
requiring that $(\tilde{u}, \tilde{U})$ satisfies~\eqref{equ:paramuoto} for all $(v,V)$ in the chosen finite-dimensional subspace of $L_P^2(\Gamma;\mathcal{H})\simeq L_P^2(\Gamma)\otimes \mathcal{H}$, i.e., for all $(v,V) = (L_{\mu'}\varphi_{j'}, L_{\mu'} \vv_{i'})$, $\mu' \in \Lambda$, $j' = 1, \dots, N_D$, $i' = 1,\dots, M-1$, one ends up at the linear system of equations (cf.~\cite{Leinonen14,Vauhkonen97})
\begin{align}
\label{equ:linearsystem}
\left(
\begin{array}{cc}
\mathbf{\Delta} & \mathbf{\Upsilon} \\
\mathbf{\Upsilon}^\mathsf{T} & \mathbf{\Pi} \\
\end{array}
\right)
\left(
\begin{array}{cc}
\alpha\\
\beta
\end{array}
\right)
= 
\left(
\begin{array}{cc}
\mathbf{0}\\
\mathbf{c}
\end{array}
\right)
.
\end{align}
Here,
$\mathbf{\Delta} \in \R^{N_D N_\Gamma \times N_D N_\Gamma}$ and $\mathbf{\Pi} \in \R^{(M-1) N_\Gamma \times (M-1) N_\Gamma}$ are symmetric sparse matrices,  
$\mathbf{\Upsilon} \in \R^{N_D N_\Gamma \times (M-1) N_\Gamma}$ is a sparse (non-square) matrix, $\mathbf{c} \in \R^{(M-1) N_\Gamma}$ is a block vector, and $\mathbf{0} \in \R^{N_D N_\Gamma}$~is a zero vector. 
Take note that~\eqref{equ:linearsystem} has in total $N_{\rm tot}:=(N_D+M-1)N_\Gamma$ degrees of freedom.

In order to give the precise definitions of the elements in the system~\eqref{equ:linearsystem}, let us first introduce some auxiliary block matrices.
In the following definitions, 
$i,i' = 1,\ldots,M-1$,
$j,j' = 1,\ldots,N_D$,
$k = 1,\ldots,L+M$,
$l = 1,\ldots,L$,
$m = 1,\ldots,M$,
and
$\mu,\mu' \in \Lambda$, if not stated otherwise.
The FEM matrices corresponding to the spatial discretization of $D$ are defined via
\begin{align*}
\{\mathbf{A}_0\}_{j,j'} &=  \int_D  \sigma_0\, \nabla \varphi_j(\x)  \cdot \nabla \varphi_{j'}(\x)  \, \dd \x \, , \\[1mm]
\{\mathbf{A}_l\}_{j,j'} &= \int_{D_l} \sigma_l\, \nabla \varphi_j(\x) \cdot \nabla  \varphi_{j'}(\x)  \, \dd \x \, .
\end{align*}
Notice that $A_0$ is sparse and $\{ A_l \}_{j,j'}$ is nonzero only if the supports of both $\varphi_j$ and $\varphi_{j'}$ intersect $D_l$. 
The elements of the stochastic moment matrices are 
\begin{align*}
\{\mathbf{G}_0\}_{\mu,\mu'} &= \int_\Gamma  L_\mu(\y)  L_{\mu'}(\y)  \, \dd \y \, = \, \delta_{\mu, \mu'},\\
\{\mathbf{G}_k\}_{\mu,\mu'} &= \int_\Gamma y_k  L_\mu(\y)  L_{\mu'}(\y) \, \dd \y.
\end{align*}
Since a univariate Legendre polynomial of a certain order is orthogonal to all lower order polynomials, it follows easily that $\{\mathbf{G}_k\}_{\mu,\mu'} \not= 0$ only if $|\mu_k - \mu'_k| = 1$ and $\mu_{k'} = \mu'_{k'}$ for $k' \not= k$, which makes $\mathbf{G}_k$ very sparse.
Finally, the electrode mass matrices are defined through
\begin{align*}
\{\mathbf{S}_m\}_{j,j'} =
\int_{E_m} \, 
\varphi_{j}(\x)\,
\varphi_{j'}(\x)\, \dd S,
\end{align*}
and the contact conductance matrices through (cf.~\eqref{sconduct})
\begin{align*}
\mathbf{Z}_m 
= \frac{1}{2}(a_m+b_m) \mathbf{G}_{0} + \frac{1}{2}(b_m-a_m) \mathbf{G}_{L+m}.
\end{align*}
Standard FEM techniques can be used to construct 
$\mathbf{A}_l$, $l=0, \dots, L$, and $\mathbf{S}_m$, $m=1, \dots, M$,
and we refer to \cite{Bieri09a,Leinonen14c}
for the efficient formation of
$\mathbf{G}_k$, $k=0, \dots, L+M$. The contact conductance matrices
$\mathbf{Z}_m$, $m=1,\dots, M$, are trivial to construct as soon as the stochastic moment matrices are available.

Now, the matrix $\mathbf{\Delta}$ can be given as
\begin{align*}
\mathbf{\Delta} = 
\sum_{l=0}^L \mathbf{A}_l \otimes \mathbf{G}_l +
\sum_{m=1}^M \mathbf{S}_m \otimes \mathbf{Z}_m,
\end{align*}
where 
$\otimes$ denotes the Kronecker product. Moreover, 
\begin{align*}
\{\mathbf{\Upsilon}_{j,i'}\}_{\mu,\mu'} = 
\{\mathbf{Z}_{i'+1}\}_{\mu,\mu'} 
\int_{E_{i'+1}}\varphi_{j}(\x) \, \dd S -
\{\mathbf{Z}_{1}\}_{\mu,\mu'}
\int_{E_1} \varphi_{j}(\x) \, \dd S
\end{align*}
and
\begin{align*}
\{\mathbf{\Pi}_{i,i'}\}_{\mu,\mu'} 
& = 
\{\mathbf{Z}_{1}\}_{\mu,\mu'}
|E_1|
+ \delta_{i,i'} \, 
\{\mathbf{Z}_{i'+1}\}_{\mu,\mu'}
|E_{i+1}|,
\end{align*}
where 
$|E_{i}|$ denotes the area/length of the $i$th electrode.
Finally, the block vector $\mathbf{c}$ is defined elementwise by
\begin{align*}
\{\mathbf{c}_{i}\}_{\mu} 
= \, 
(I \cdot \vv_{i}) \,
\int_\Gamma L_{\mu}(\y) \dd \y
 &= \begin{cases}
   0, &  \mu \neq \mathbf{0}, \\
   I_1-I_{i+1}, & \mu = \mathbf{0},
  \end{cases}
\end{align*}
where $I \in \R^M_{\diamond}$ is the applied current pattern and $\mathbf{0}$ is the zero multi-index.

\section{Inverse solution}
\label{sec:inverse_solution}

The objective of EIT is to retrieve useful information about the conductivity inside the examined body based on measured noisy electrode current-potential pairs.
In this section, we explain how the sGFEM approximation \eqref{sGFEM_U} for the second component of the solution to \eqref{equ:paramuoto} can be employed in numerically solving this problem in the Bayesian framework; see \cite{Kaipio04a} for more information on statistical inversion. 

Let $I^1, \dots, I^{M-1} \in \R^M_\diamond$ be linearly independent current patterns that are driven in turns through the $M$ contact electrodes $E_1, \dots, E_M$, and suppose $V^1, \dots, V^{M-1} \in \R^M$ are the corresponding measured noisy electrode potential vectors. (Notice that there is no benefit in using more than $M-1 = {\rm dim}(\R^M_\diamond)$ current patterns because the solution of \eqref{equ:stokamuoto} depends linearly on $I$.) We define 
$$
\vbf = \Big[(V^1)^{\mathsf{T}}, \dots, (V^{M-1})^{\mathsf{T}}\Big]^{\mathsf{T}} \in \R^{M(M-1)}
$$ 
and 
$$
\tilde{\mathcal{U}}(\mathbf{Y}) = \Big[\tilde{U}^1(\mathbf{Y})^{\mathsf{T}}, \dots, \tilde{U}^{M-1}(\mathbf{Y})^{\mathsf{T}}\Big]^{\mathsf{T}} \in \R^{M(M-1)}
$$ 
with $\tilde{U}^m(\mathbf{Y}) \in \R^M_\diamond$ being the sGFEM solution \eqref{sGFEM_U} corresponding to the current pattern $I = I^m$ in \eqref{equ:paramuoto}.
In other words, $\tilde U^i_j(\mathbf{Y}) \in \R$
is the $j$th component of the sGFEM solution \eqref{sGFEM_U} for the current pattern $I^i \in \R^M_\diamond$.

The electrode potentials $\vbf$ are assumed to be a realization of the random variable
\begin{align}
\label{equ:Bayesmodel}
\mathbf{V} = \tilde{\mathcal{U}}(\mathbf{Y}) + \mathbf{E},
\end{align}
where $\mathbf{E}$ is the noise process contaminating the measurements. Notice that the model \eqref{equ:Bayesmodel} cannot be exact as it does not take into account the unavoidable discretization errors in $\tilde{\mathcal{U}}(\mathbf{Y})$, but we choose to ignore this fact to simplify the analysis.
Moreover, $\mathbf{E}: \Omega \to \R^{M(M-1)}$ is assumed to be independent of $\mathbf{Y}$, mean-free, and Gaussian with a known covariance matrix $\mathbf{L} \in \R^{M(M-1) \times M(M-1)}$.
Combining \eqref{equ:Bayesmodel} with the Bayes' formula results in the posterior density
\begin{align}
\label{Bayes1}
\pi(\y \, | \, \vbf) \, &\propto  \, \pi_{\rm noise}\big(\vbf - \tilde{\mathcal{U}}(\mathbf{y}) \big) \, \pi_{\rm pr}(\y) \nonumber\\[1mm]
& = \frac{1}{\sqrt{(2\pi)^{M(M-1)}|\mathbf{L}|}}\exp\!\Big(\!-\frac{1}{2}(\vbf - \tilde{\mathcal{U}}(\mathbf{y}))^{\transposeT}\mathbf{L}^{-1}(\vbf - \tilde{\mathcal{U}}(\mathbf{y}))\Big) \, \pi_{\rm pr}(\y),
\end{align}
where $|\mathbf{L}|$ is the determinant of the noise covariance matrix and the `constant' of proportionality is independent of $\y$.

The choice of the prior density $\pi_{\rm pr}$ in \eqref{Bayes1} should be based on {\em a priori} information about the pixel values of the conductivity and the contact conductances. Since the sGFEM forward solver of the previous section was already built under the assumption that the parameters $\y$ belong to the hypercube $\Gamma = [-1,1]^{L+M}$, it is natural to choose
\begin{equation}
\label{prior}
\pi_{\rm pr}(\y) \, = \, \pi_\sigma(\y_{\sigma}) \pi_\zeta(\y_{\zeta}) \mathbf{1}_\Gamma(\y),
\end{equation}
where $\mathbf{1}_\Gamma: \R^{L+M} \to \R$ is the indicator function of $\Gamma \subset \R^{L+M}$ and we have assumed that the parameters corresponding to the pixelwise conductivity values $\y_{\sigma} \in \R^L$ and those associated to the contact conductances $\y_\zeta \in \R^M$ are independent {\em a priori}. We assume to have no further prior information on the contact conductances, i.e., we employ
$$
\pi_\zeta(\y_\zeta) \, = \, 2^{-M}, \qquad  \y_\zeta \in  [-1,1]^M,
$$
whereas for the conductivity we choose a truncated multivariate normal prior density:
\begin{align}
\label{prior2}
\pi_{\sigma}(\y_\sigma) = 
\frac
{ \exp \Big( - {\displaystyle \frac{1}{2}} {\y_\sigma^\transposeT} \mathbf{M}^{-1} \y_{\sigma} \Big) }
{ \displaystyle{\int_{\Gamma_\sigma}} \exp \Big( - \frac{1}{2} \tilde\y_\sigma^\transposeT\,\mathbf{M}^{-1}\,\tilde\y_{\sigma} \Big) \dd \tilde\y_\sigma },  \qquad
\y_\sigma \in  [-1,1]^L,
\end{align}
where 
$\Gamma_\sigma = [-1,1]^L$ and 
$\mathbf{M}\in\R^{L \times L}$ is the covariance matrix of the underlying multivariate normal distribution $\mathcal{N}(\mathbf{0},\mathbf{M})$.
In this work, the covariance matrix is assumed to be of the squared exponential type:
\begin{align}
\label{equ:priorcov}
\mathbf{M}_{l,l'} = \eta^2 \exp \left( \frac{-|\mathbf{r}_l-\mathbf{r}_{l'}|^2}{2s^2} \right),
\end{align}
where $\mathbf{r}_l$ is the center of the pixel $D_l$, $s>0$ is the correlation length, $\eta > 0$ is the standard deviation, and $l,l'=1,\ldots,L$.

\begin{remark}
The inclusion of $\mathbf{1}_\Gamma(\mathbf{y})$ in \eqref{prior} is only natural because there is {\em absolutely} no guarantee that $\tilde{\mathcal{U}}(\mathbf{y})$ is any kind of an approximation for the electrode potentials corresponding to a conductivity of the form \eqref{equ:klexpansion} if $\mathbf{y} \notin \Gamma$. The `additional' priors $\pi_\sigma$ and $\pi_\zeta$ can, however, be selected as one wishes,
bearing in mind that complicated choices may hamper the computation of the MAP and CM estimates for the posterior.

One could
also utilize the prior information in $\pi_\sigma$ and $\pi_\zeta$ when building the sGFEM forward solver to maximize the accuracy of $\tilde{\mathcal{U}}(\mathbf{y})$ for those parameter vectors $\mathbf{y}$ that live in regions of high prior probability (cf.~\cite{Hakula14}).
One way of achieving this is to replace the probability density \eqref{equ:probden} by an approximation of \eqref{prior} in the Legendre polynomial basis and use techniques in \cite{Leinonen14c} to construct the (more involved) stochastic moment matrices.

 The reason for not taking such a path in this work is two-fold: (i) Changing $\pi_\sigma$ and $\pi_\zeta$ does not affect the sGFEM forward solver in our setting, which significantly reduces the computational cost for tuning/changing the prior. (ii) Using a more complicated random field model than \eqref{equ:klexpansion0} for the sGFEM forward solver leads easily to a less sparse system matrix \eqref{equ:linearsystem} that is more laborious to construct, and it potentially also makes controlling the positivity of the conductivity more involved. 

\end{remark}

The MAP estimate $\y_{\rm MAP}$ for $\mathbf{Y}$, i.e., the maximizer of the posterior density \eqref{Bayes1}, can be computed by solving the constrained minimization problem
\begin{align}
\label{equ:MAPproblem}
\y_{\rm MAP} \, := \, \operatorname*{arg\,min}_{\y\in\Gamma} F(\y),
\end{align}
where 
\begin{align*}
F(\y) := \big(\vbf - \tilde{\mathcal{U}}(\y)\big)^{\transposeT} \mathbf{L}^{-1} \big(\vbf - \tilde{\mathcal{U}}(\y)\big) + \y_\sigma^\transposeT\,\mathbf{M}^{-1}\,\y_{\sigma}
\end{align*}
is a positive-valued polynomial in $\y$.
Subsequently, the MAP estimate for the conductivity $\sigma_{\rm MAP}: D \to \R_+$ is obtained by evaluating \eqref{equ:klexpansion} at $\y = \y_{\rm MAP}$, and the MAP estimates for the contact conductances are deduced analogously via \eqref{sconduct}.

The CM estimates of the conductivity and contact conductances are obtained by (numerically) evaluating the $(L+M)$-dimensional integrals
\begin{align}
\label{CMint1}
\sigma_{\rm CM}(\x) = \int_{\Gamma} \sigma(\y,\x) \pi(\y \, | \, \vbf) \dd \y, \qquad \x \in D,
\end{align}
and
\begin{align}
\label{CMint3}
(\zeta_m)_{\rm CM} = \int_{\Gamma} \zeta_m(\y) \pi(\y \, | \, \vbf) \dd \y, \qquad m = 1,\ldots,M,
\end{align}
respectively.
To evaluate the reliability of the CM estimates, we also consider the conditional {\em standard deviations} (SD)
\begin{align}
\label{CMint5}
\sigma_{\rm SD}(\x) = \sum_{l=1}^L \sigma_l \mathbf{1}_{D_l}(\x) \left[ \int_\Gamma y_l^2\, \pi(\y \, | \, \vbf) \dd \y - \left(\int_\Gamma y_l \,  \pi(\y \, | \, \vbf) \dd \y\right)^2 \right]^{\frac{1}{2}}, \qquad \x \in D,
\end{align}
and
\begin{align}
\label{CMint6}
(\zeta_m)_{\rm SD} = \frac{1}{2}(b_m-a_m)\left[ \int_\Gamma y_{L+m}^2\, \pi(\y \, | \, \vbf) \dd \y - \left(\int_\Gamma y_{L+m} \,  \pi(\y \, | \, \vbf) \dd \y\right)^2 \right]^{\frac{1}{2}}
\end{align}
in the numerical experiments of Section~\ref{sec:numerical_examples}.

\section{Two-phase implementation}
\label{sec:implementation}
The implementation of the presented inversion algorithm consists of two phases:
the pre-measurement and post-measurement processing.
The former corresponds to computations that can be carried out before performing any measurements, assuming the object shape, the electrode positions, and the preliminary bounds for the conductivity and contact conductances are known. The latter phase consists of forming the posterior density and computing the desired estimates for the unknowns.

\subsection{Pre-measurement processing}

The pre-measurement phase consists of the following six steps:
\begin{enumerate}
\item Specify the computational domain, i.e., the object shape together with the electrode sizes and positions.
\item Select a suitable partition of the domain into pixels.
\item Specify bounds for the conductivity and contact conductance values, i.e., $\sigma_0, \sigma_1, \dots, \sigma_L$ in \eqref{equ:klexpansion} as well as $a_1, \dots , a_M$ and $b_1, \dots , b_M$ in \eqref{sconduct}.
\item Construct a suitable FEM polynomial basis for $H^1(D)$.
\item Select the index set $\Lambda$ for the polynomial chaos expansion.
\item Compute the sGFEM solution \eqref{sGFEM_u}--\eqref{sGFEM_U}. 
\end{enumerate}
We emphasize that all these steps can be performed without having the actual electrode measurements in hand. Moreover, the sGFEM solution can be reused for different data sets as long as the bounds for the conductivity and contact conductances or the measurement geometry are not altered.

The pre-measurement processing stage is clearly the more time consuming of the two phases because the SCEM forward problem is discretized by over $10^7$ degrees of freedom in our two-dimensional numerical experiments. 
(In three dimensions, the number of degrees of freedom could easily exceed $10^9$.) Fortunately, if the measurement configuration is known well in advance, the pre-measurement processing can be carried out before the actual measurements.

\subsection{Post-measurement processing}
After the electrode potential measurements $\vbf \in \R^{M(M-1)}$ are available, the post-measurement phase consists of the following four steps:
\begin{enumerate}
\item Specify the noise covariance matrix $\mathbf{L}$.
\item Select the correlation length $s$ and the standard deviation $\eta$ for the prior covariance matrix $\mathbf{M}$ in \eqref{equ:priorcov}.
\item Construct the posterior density \eqref{Bayes1}.
\item Compute the desired estimates (MAP, CM, and SD) for the posterior distribution.
\end{enumerate}

Notice that the accuracy of the spatial FEM discretization does not affect the computation time for the post-measurement phase since the approximate stochastic forward solution $\tilde{U}(\y)$ from \eqref{sGFEM_U} does not involve the spatial FEM basis functions. Hence, one should use as dense spatial FEM mesh as allowed by the {\em pre-measurement} time and memory constraints. 
On the other hand, the discretization of $L^2(\Gamma)$ affects the computation times of both phases.

\section{Numerical experiments}
\label{sec:numerical_examples}

\begin{figure}[t]
  \begin{center}
  {\includegraphics[]{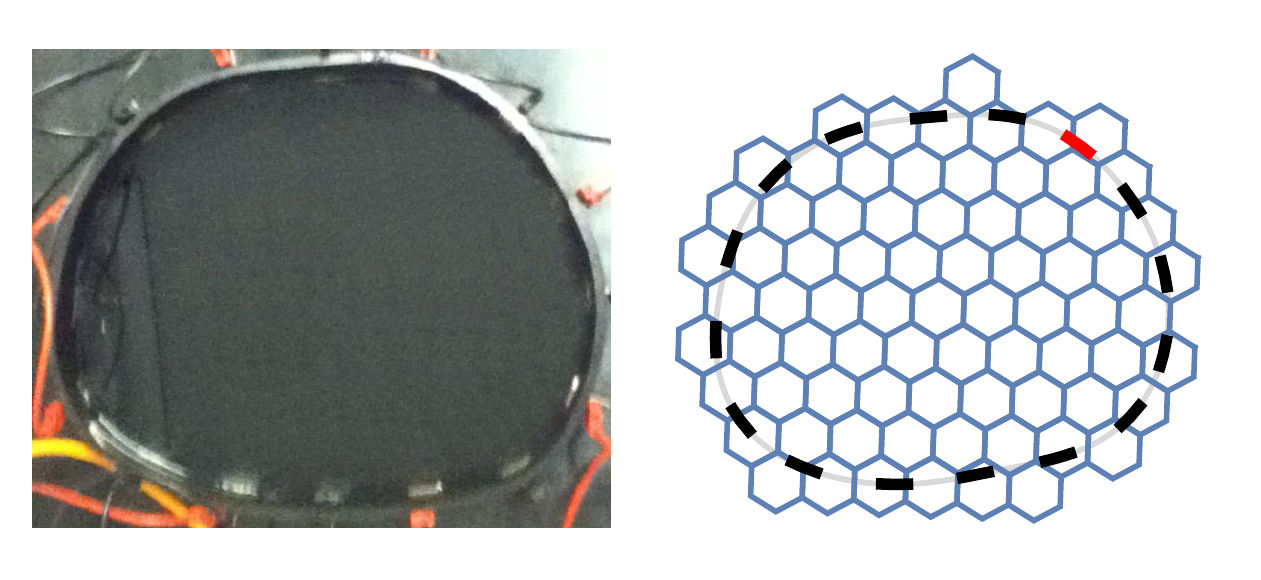}}
  \end{center}
  \caption{
Left: thorax-shaped water tank with no inclusions.
Right: the computational domain, its partition into $L=76$ hexagonal pixels, and the $M=16$ attached electrodes. The current-feeding electrode $E_1$ is red and the others are numbered in counterclockwise order.
}
  \label{fig:pixelgrid}
\end{figure}

We apply the above introduced methodology to five sets of experimental data from a thorax-shaped water tank with vertically homogeneous embedded objects of steel and/or plastic extending from the bottom all the way through the water surface.  The circumference of the tank is $106\,{\rm cm}$, and $M=16$ rectangular metallic electrodes of width $2\,{\rm cm}$ and height $5\,{\rm cm}$ are attached to the interior lateral surface of the tank. In all tests, the tank is filled with tap water up to the top of the electrodes.  The measurement configuration without inclusions is presented in the left-hand image of Figure~\ref{fig:pixelgrid}.
(All photographs shown below are cropped and spatially normalized versions of the original ones. We have also removed most of the reflections on the water surface to ease perceiving the images.)
The measurements were performed with low-frequency ($1\,{\rm kHz}$) alternating current using the {\em Kuopio impedance tomography} (KIT4) device~\cite{Kourunen09}. The phase information of the measurements is ignored, meaning that the amplitudes of electrode currents and potentials are interpreted as real numbers. The employed (real) current patterns are (cf.~\eqref{meanfreebasis}) 
$$
I^m =  ({\mathrm e}_1 - {\mathrm e}_{m+1}) \, {\rm mA} , \qquad m = 1, \dots, M-1,
$$
with ${\mathrm e}_m$ denoting the $m$th Euclidean basis vector of $\R^M$. This choice of current basis makes the first electrode special; it is marked with red color in Figure~\ref{fig:pixelgrid}.

As the measurement setting is vertically homogeneous --- notice that no current flows through the bottom or the top of the water tank, which corresponds to homogeneous Neumann boundary conditions --- it can be modeled with a two-dimensional version of the SCEM (cf.~Section~\ref{sec:SCEM}).
The conversion of conductivity (${\rm S/m}$) and contact conductances (${\rm S/m}^2$) into corresponding two-dimensional quantities is achieved by multiplying with the height of the electrodes.
The same measurement setting was tackled in \cite{Darde13b}, where the conductivity of tap water was estimated to be around $0.2$ -- $0.25\,{\rm mS/cm}$, i.e., $1.0$ -- $1.25\,{\rm mS}$ in the two-dimensional units. This also matches the limits given for drinking water in the literature ($0.05$ -- $0.5\,{\rm mS/cm}$).
Using \cite{Darde13b} as our reference, we choose $\sigma_0 = 1.1\, {\rm mS}$ and $\sigma_1, \dots  , \sigma_L = 0.9\, {\rm mS}$ in \eqref{equ:klexpansion}, i.e., we let the pixelwise conductivities vary between $0.2\,{\rm mS}$ and $2.0\,{\rm mS}$ in the forward solver. As the examples consider inclusions that are either insulating (plastic) or highly conducting (steel), the interval $[0.2, 2.0] \,{\rm mS}$ for the conductivity values may seem a bit restrictive. However, according to our experience (cf.,~e.g.,~\cite{Darde13b,Harhanen15}), $0.2\,{\rm mS}$ is a sufficiently low value for modeling an insulating object accurately enough and, on the other hand, highly conducting objects exhibit some resistivity in EIT, probably due to the contact resistance at their boundaries (cf.~\cite{Heikkinen01}). A relatively large lower bound for the conductivity also ensures that the sGFEM system matrix stays well conditioned. 
Furthermore, we assume relatively bad contacts at the electrode-water interfaces and set $a_m = 10\, {\rm mS/cm}$ and $b_m = 10^3\, {\rm mS/cm}$, $m = 1,\ldots,M$ in \eqref{sconduct} (cf.~\cite{Heikkinen02}).

The right-hand image of Figure~\ref{fig:pixelgrid} shows the computational domain $D\subset \R^2$ corresponding to the water tank together with our choice for the partition of the domain into $L=76$ pixels $D_1, \dots , D_L$ (cf.~\eqref{equ:klexpansion}) that are intersections of certain hexagons and $D$. We employ spatial FEM mesh (not shown) composed of $N_D=9383$ nodes with appropriate refinements at the edges of the electrodes (cf.~\cite{Darde13b}). As the stochastic index set in \eqref{sGFEM_u}--\eqref{sGFEM_U}, we use $\Lambda = \isoTD(L + M, 2)$, which results in $N_\Gamma = 4371$ stochastic degrees of freedom.
In total, the discretized forward SCEM problem includes $N_{\rm tot} = (N_D+M-1)N_\Gamma \approx 4.1 \cdot 10^7$ unknowns, and the system matrix in \eqref{equ:linearsystem} has approximately $3\cdot 10^8$ nonzero elements, i.e., approximately seven nonzero elements per row. In all our numerical experiments, \eqref{equ:linearsystem} is solved by the standard direct linear solver of MATLAB,
i.e., by the {\tt mldivide} command, for simplicity and to avoid any convergence and preconditioning issues related to iterative methods. Using the conjugate gradient method with an ILU0 \cite{Saad03} based preconditioner, we have been able to tackle denser FEM and pixel meshes, e.g., $L=145$ corresponding to $N_\Gamma = 13203$ and $N_{\rm tot} \approx 1.2 \cdot 10^8$, 
but this does not result in significantly better results than the ones presented in Sections~\ref{sec:empty}--\ref{sec:two} below.

To motivate the choice of the stochastic index set, we mention that for $\Lambda = \isoTD(L+M,1)$, the conductivity reconstructions contain more artifacts, the inclusions are not as well localized, and the background conductivity level is higher and not as smooth as with $\isoTD(L+M,2)$. We were not able to test the case $\Lambda = \isoTD(L+M,3)$ with any reasonable FEM and pixel meshes due to memory and time constraints. 
There is an obvious trade-off between the fineness of the FEM mesh and the number of the hexagonal pixels in the reconstruction grid; the values listed in the previous paragraph represent a compromise arrived at via trial and error.
Employing denser FEM mesh forces one to use a coarser pixel grid --- and vice versa --- in order to keep the system size reasonable.
Take note that increasing the number of spatial degrees of freedom $N_D$ affects the size of the sGFEM system \eqref{equ:linearsystem} linearly, whereas increasing $L$ leads to a quadratic growth rate since
$$
N_\Gamma = {2+(L+M) \choose 2} = \frac{(L+17)(L+18)}{2}
$$
for $\Lambda = \isoTD(L+M,2)$ and $M=16$ electrodes (cf.~\eqref{stoch_df}). Recall also that increasing $N_D$ affects only the computation time of the pre-measurement stage while the number of pixels in the reconstruction grid has an effect on the time consumption in both pre- and post-measurement phases.

The magnitude of the measurement noise on each electrode is assumed to be proportional to the difference of the smallest and largest electrode potential measurement, leading to the choice (cf.~\cite{Darde13b}) 
\begin{equation}
\label{noisevar}
\mathbf{L} = \xi^2 \mathbf{I},\,
\qquad
\xi = 0.01 \, (\max(\vbf) - \min(\vbf)),
\end{equation}
for the noise covariance matrix. Here and in the following, $\mathbf{I}$ denotes an identity matrix of the appropriate size. Loosely speaking, \eqref{noisevar} corresponds to assuming one per cent of measurement noise. As the noise level of the measurement device is probably only a couple of per mille depending on the measurement channel \cite{Kourunen09}, the assumed high variance for the noise process is actually used partially to mask the unavoidable discretization errors in the sGFEM forward solution for \eqref{equ:stokamuoto}; see \cite[Remark~5.1]{Hakula14}.  
We use the correlation length $s = 5$ and the standard deviation $\eta = 10\xi$ in the prior covariance matrix $\mathbf{M}$ of \eqref{equ:priorcov}. The choice of $s$ reflects the prior assumption on the diameter of the embedded inhomogeneities, while the values of the other free parameters $\xi$ and $\eta$ were chosen by trial and error, guided by the last test case (cf.~Figure~\ref{fig:Example1_5}).
The prior covariance matrix was constructed assuming that all pixels are hexagonal, and hence some center points of the pixels actually lie outside the computational domain.

The MAP estimate $\y_{\rm MAP}$ --- and subsequently $\sigma_{\rm MAP}$ --- is obtained by solving \eqref{equ:MAPproblem} as a nonlinear least-squares minimization problem by resorting to
the {\tt lsqnonlin} function provided by the Optimization Toolbox of MATLAB. 
The CM estimates for the conductivity and the contact conductances as well as the related standard deviations are computed via Markov chain Monte Carlo (MCMC) simulations; the usage of a deterministic sparse quadrature rule such as the one of Smolyak \cite{Smolyak63,Schillings13} would be another possibility, but we have had more success with MCMC techniques in connection with EIT.
The standard Metropolis--Hastings algorithm (see,~e.g.,~\cite{Kaipio04a}) is used to generate a sample of parameter vectors
\begin{align*}
\{ \y^{(1)},\ldots,\y^{(N)}\} \subset \R^{L+M}
\end{align*}
that is distributed (approximately) according to the posterior $\pi(\y~|~\vbf)$ given by~\eqref{Bayes1}.
Starting from the corresponding MAP estimate, we use a single random walk, with a burn-in period of $5 \cdot 10^4$ and a
thinning of five, i.e., we only store every fifth element of the Markov chain, to generate $N=4 \cdot 10^5$ samples. The proposal density for the random walk is 
the truncated multivariate normal on $\Gamma$
centered at the previous sample with the covariance matrix $0.07^2\, \mathbf{I}$, resulting in an acceptance rate of approximately $30\%$.
Subsequently, the integrals \eqref{CMint1} and \eqref{CMint3} are approximated as 
\begin{align*}
\sigma_{\rm CM}(\x) \approx \frac{1}{N} \sum_{i=1}^N \sigma(\y^{(i)},\x) \quad \textrm{ and } \quad 
(\zeta _m)_{\rm CM} \approx \frac{1}{N} \sum_{i=1}^N \zeta_m(\y^{(i)}),
\end{align*}
respectively.
Similarly, the standard deviations \eqref{CMint5} and \eqref{CMint6} are approximated as
\begin{align*}
\sigma_{\rm SD}(\x) \approx \sum_{l=1}^L \sigma_l \mathbf{1}_{D_l}(\x)
\left[ 
\frac{1}{N} \sum_{i=1}^N (y^{(i)}_l)^2 - \left(\frac{1}{N} \sum_{i=1}^N y^{(i)}_l \right)^2
\right]^{\frac{1}{2}}
\end{align*}
and
\begin{align*}
(\zeta_m)_{\rm SD} \approx \frac{1}{2}(b_m-a_m)
\left[
\frac{1}{N} \sum_{i=1}^N (y^{(i)}_{L+m})^2 - \left(\frac{1}{N} \sum_{i=1}^N y^{(i)}_{L+m} \right)^2
\right]^{\frac{1}{2}},
\end{align*}
respectively. The number of samples was evaluated to be sufficient by visually examining the development of the CM estimates: in all numerical examples, the estimates seemed to stabilize after about $2 \cdot 10^5$ samples --- the final sample size was chosen to be twice as large.

The solution of the SCEM forward problem and most other computations were performed using the commercial software packages MATLAB\footnote{Version 8.2.0 (R2013b), The MathWorks Inc., Natick, Massachusetts, 2013.} 
and Mathematica\footnote{Version 9.0, Wolfram Research Inc., Champaign IL, 2012.}.
MATLink~\cite{MATLink14} 
was employed for seamless two-way communication and data transfer between Mathematica and MATLAB, and the needed FEM meshes were generated by NETGEN mesh generator~\cite{NETGEN10}.

\begin{figure}[t!]
  \begin{center}
  {\includegraphics[]{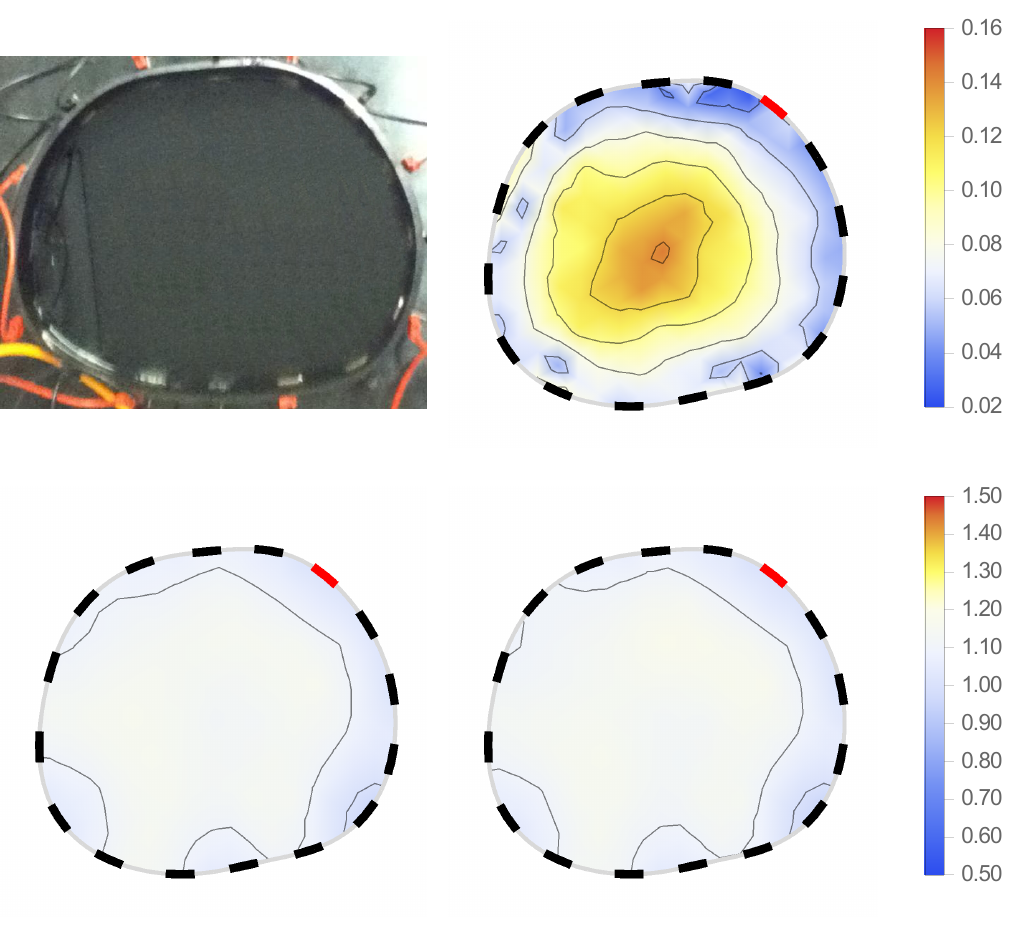}}
  \end{center}
  \caption{
    Results of the first example.
    Top left: the target without embedded inclusions.
    Top right: the SD estimate $\sigma_{\rm SD}$.
    Bottom left: the MAP estimate $\sigma_{\rm MAP}$.
    Bottom right: the CM estimate $\sigma_{\rm CM}$.
    The unit in all images is mS, the MAP and CM estimates use the same colormap, and the colorbar tick markers correspond to the contour lines in the images.
}
  \label{fig:Example1_0}
\end{figure}

\subsection{Experiment with empty tank}
\label{sec:empty}

As a first simple example, we consider the setting in the top left image of Figure~\ref{fig:Example1_0}, i.e., the case of no embedded inclusions. The other images of Figure~\ref{fig:Example1_0} show interpolated versions of the pixelwise SD, MAP, and CM estimates for the conductivity. 
Both MAP and CM estimates produce tolerable and almost identical reconstructions of the empty tank. Take note that some of the small artifacts close to the object boundary are probably caused by mismodeled geometry: the shape of the water tank and the positions of the electrodes  were estimated based on the photographs and previous experiments with the same measurement configuration (cf.~\cite{Darde13b}).
As expected, the SD estimate reveals that the degree of uncertainty in the conductivity reconstruction is the highest in the central parts of the tank and the lowest by the object boundary, with the smallest values of $\sigma_{\rm SD}$ occurring close to the current-feeding (red) electrode. 
The SD estimates in the other four test cases follow this same intuitive pattern.

\begin{table}[t!]
\footnotesize
\begin{center}

\begin{tabular}{r c c c c c c c c c c c c c c c c c c c c}
\\Electrode&\bf{1}&\bf{2}&\bf{3}&\bf{4}&\bf{5}&\bf{6}&\bf{7}&\bf{8} \\
\hline
MAP&10&10&17&10&202&137&800&11\\
CM&448&546&567&564&340&497&446&453\\
SD&288&279&275&270&251&253&237&291\\
~\\
Electrode&\bf{9}&\bf{10}&\bf{11}&\bf{12}&\bf{13}&\bf{14}&\bf{15}&\bf{16}\\
\hline
MAP&382&681&11&13&10&27&10&762\\
CM&587&430&654&479&523&309&513&577\\
SD&241&250&268&265&283&232&289&282
\end{tabular}
\end{center}
\caption{
The MAP, CM, and SD estimates for the contact conductances in the first experiment (mS/cm).
The mean values of these MAP, CM, and SD estimates over the sixteen electrodes are 193, 496, and 266 mS/cm, respectively.
}
\label{tab:contact_conductances_ex1}
\end{table}

\begin{figure}[t!]
  \begin{center}
  {\includegraphics[]{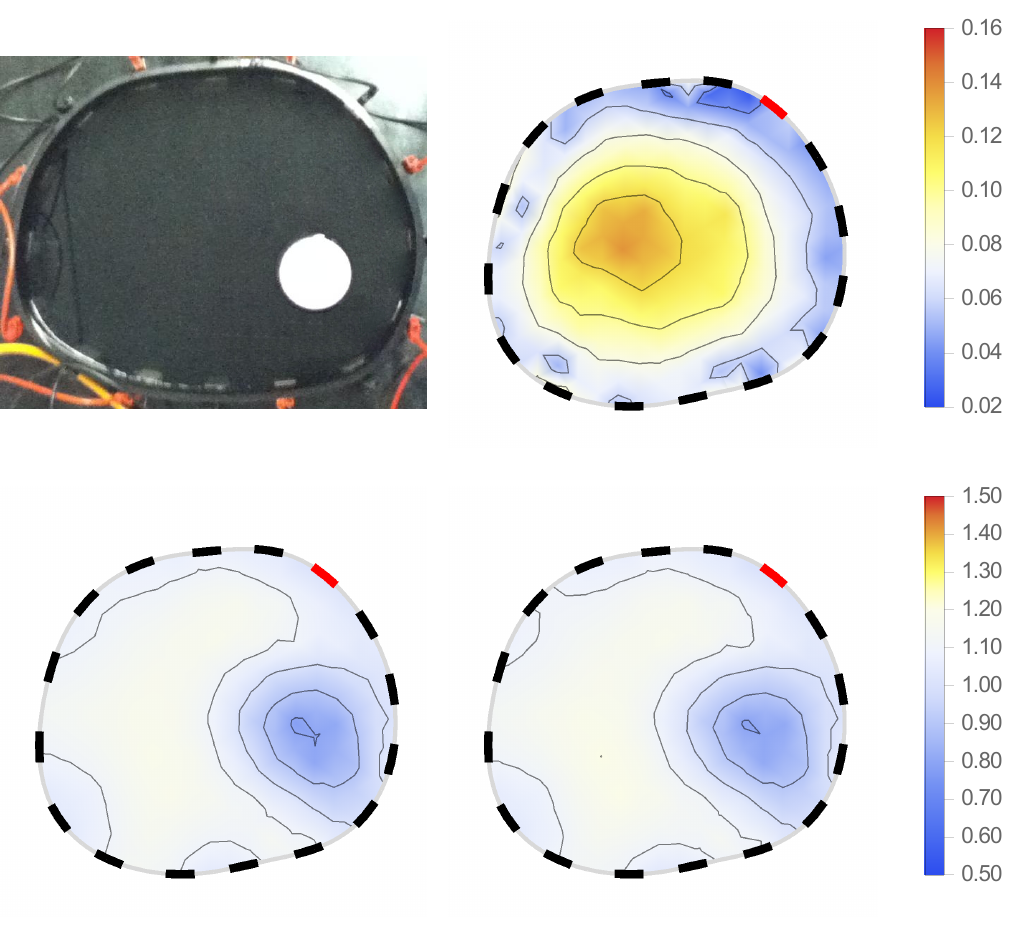}}
  \end{center}
  \caption{
    Results of the second example.
    Top left: the target with one embedded insulating inclusion.
    Top right: the SD estimate $\sigma_{\rm SD}$.
    Bottom left: the MAP estimate $\sigma_{\rm MAP}$.
    Bottom right: the CM estimate $\sigma_{\rm CM}$.
    The unit in all images is mS, the MAP and CM estimates use the same colormap, and the colorbar tick markers correspond to the contour lines in the images.
}
  \label{fig:Example1_1}
\end{figure}

The contact conductance estimates for the first experiment are presented in Table~\ref{tab:contact_conductances_ex1}.
For most electrodes, the MAP estimates of the contact conductances are close to the allowed minimum value, whereas the CM estimates stay at a higher level. One possible explanation for the low MAP estimates is the algorithm's attempt to explain the overall resistivity of the tank by introducing as high contact resistances as possible --- recall that we introduced no additional prior for the contact conductances in the post-measurement phase.
Both the MAP and CM estimates give mean contact conductances that are below the center of the interval $[10, 10^3]\, {\rm mS/cm}$ assumed in the sGFEM forward solver; see Table~\ref{tab:contact_conductances_ex1}.
We do not consider contact conductance estimates in the remaining examples as the general conclusions are the same as in this preliminary test --- and because the estimates for the contact conductances are not as interesting as the reconstructions of the conductivity.

\subsection{Experiments with one inclusion}
\label{sec:one}

\begin{figure}[t!]
  \begin{center}
  {\includegraphics[]{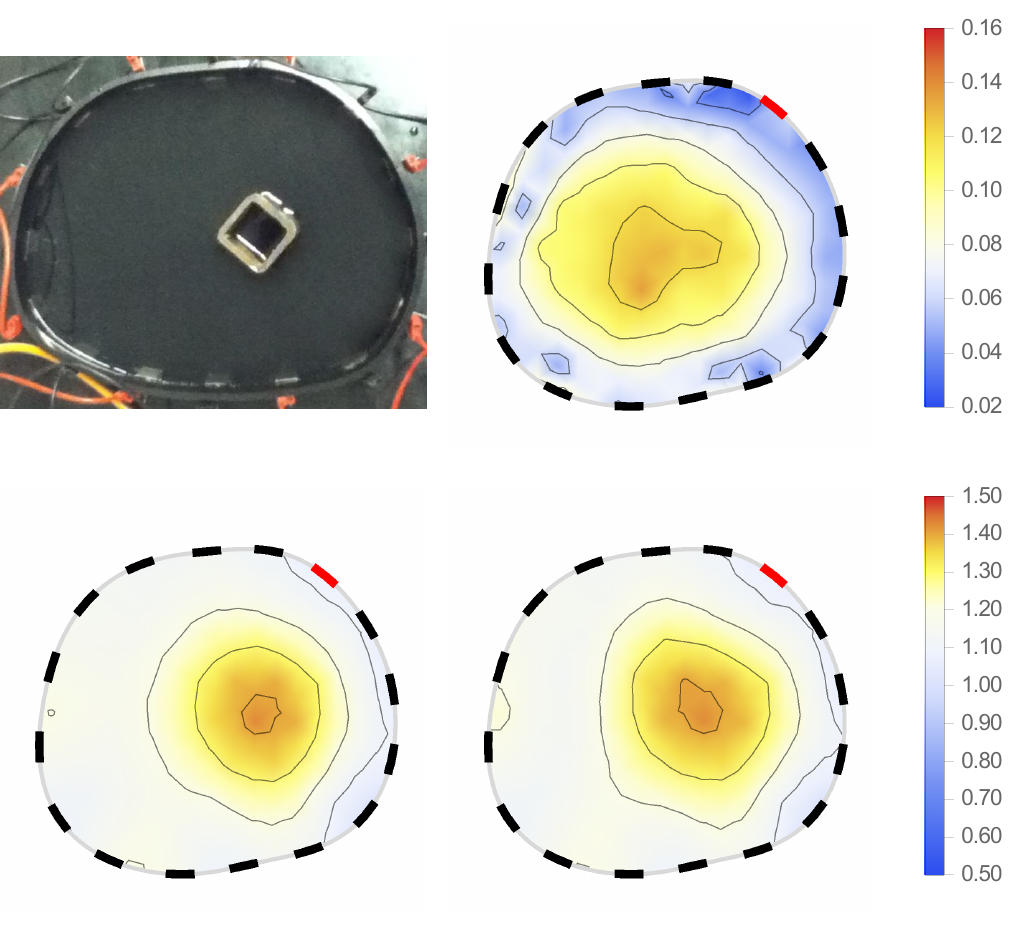}}
  \end{center}
  \caption{
    Results of the third example.
    Top left: the target with one embedded highly conducting inclusion.
    Top right: the SD estimate $\sigma_{\rm SD}$.
    Bottom left: the MAP estimate $\sigma_{\rm MAP}$.
    Bottom right: the CM estimate $\sigma_{\rm CM}$.
    The unit in all images is mS, the MAP and CM estimates use the same colormap, and the colorbar tick markers correspond to the contour lines in the images. 
}
  \label{fig:Example1_3}
\end{figure}

The top left image of Figure~\ref{fig:Example1_1} shows the target configuration of the second experiment: one insulating plastic cylinder embedded in the bottom right corner of the water tank.
The other images in Figure~\ref{fig:Example1_1} are organized as in Figure~\ref{fig:Example1_0}, and they portray the MAP, CM, and SD estimates for the conductivity. Both the MAP and CM estimates are able to find the general location of the cylinder, with the MAP estimate providing a slightly better localization. In the third experiment, one hollow steel cylinder with rectangular cross-section is immersed in the water tank; see the top left image of Figure~\ref{fig:Example1_3}. The MAP and CM estimates presented in the bottom row of Figure~\ref{fig:Example1_3} provide reasonable reconstructions of the phantom also in this case, with the hump in the MAP estimate being once again slightly sharper than in the CM estimate.
Notice that the minimal and maximal conductivity levels in the MAP and CM estimates of Figures~\ref{fig:Example1_1} and \ref{fig:Example1_3} do not lie close to the respective end points of the pixelwise interval $[0.2, 2.0]\, {\rm mS}$ used in the sGFEM forward solver: the Gaussian smoothness prior \eqref{prior2} employed in the post-measurement phase of the algorithm considerably restricts the spatial variations in the reconstructions of the conductivity.

A comparison of the reconstructions in Figures~\ref{fig:Example1_1} and \ref{fig:Example1_3} reveals that inclusions close to the exterior boundary are better localized than those deep inside the domain,
which is not surprising taking into account the general form of the SD estimates.  This trend does not depend significantly on the type of the inclusion (insulating or highly conducting) or its location in relation to the current-feeding electrode.
Notice that the correlation length $s=5\, {\rm cm}$ in the prior covariance matrix \eqref{equ:priorcov} is arguably somewhat conservative: we also tested smaller values such as $s = 3\, {\rm cm}$, which typically resulted in better resolution and contrast for the (target) inclusions, but in some cases small inclusion-like artifacts also appeared in the background,~i.e.,~at locations where there is only water inside the tank.

\begin{figure}[t!]
  \begin{center}
  {\includegraphics[]{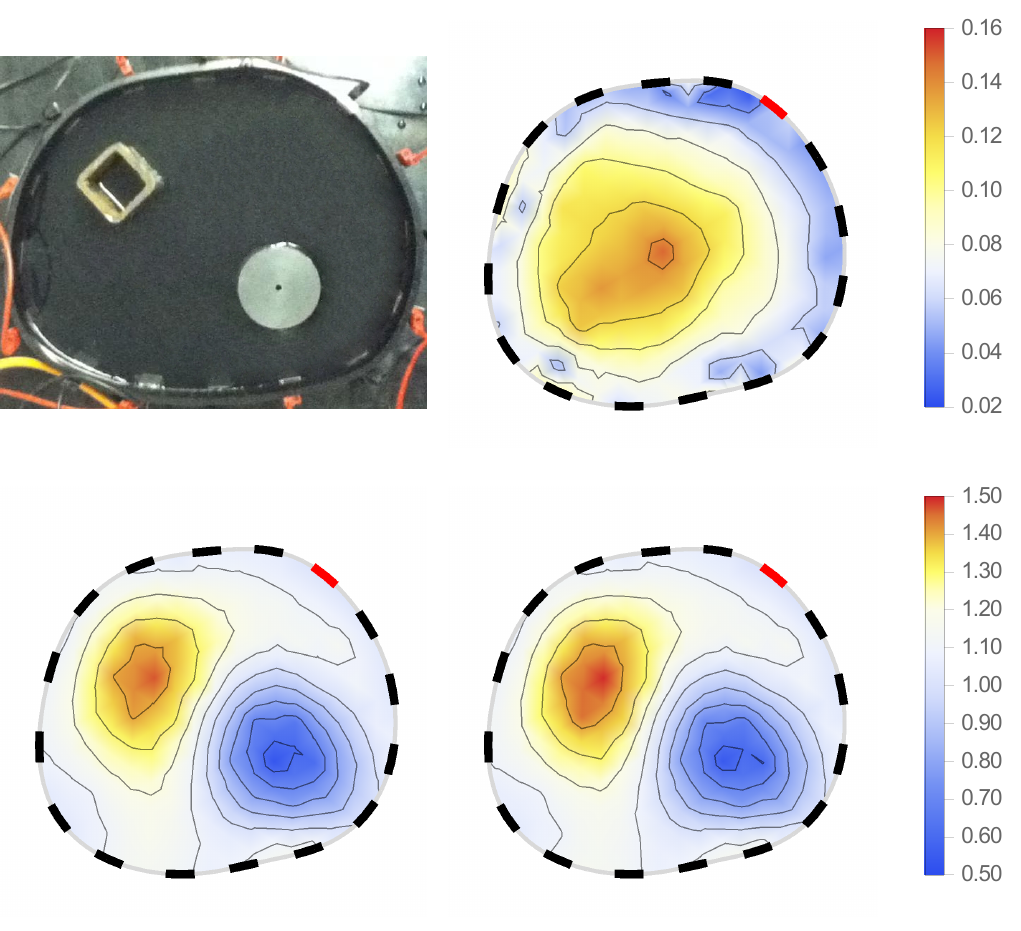}}
  \end{center}
  \caption{
    Results of the fourth example.
    Top left: the target with one insulating and one highly conducting inclusion.
    Top right: the SD estimate $\sigma_{\rm SD}$.
    Bottom left: the MAP estimate $\sigma_{\rm MAP}$.
    Bottom right: the CM estimate $\sigma_{\rm CM}$.
    The unit in all images is mS, the MAP and CM estimates use the same colormap, and the colorbar tick markers correspond to the contour lines in the images. 
}
  \label{fig:Example1_4}
\end{figure}

\subsection{Experiments with two inclusions}
\label{sec:two}
\begin{figure}[t!]
  \begin{center}
  {\includegraphics[]{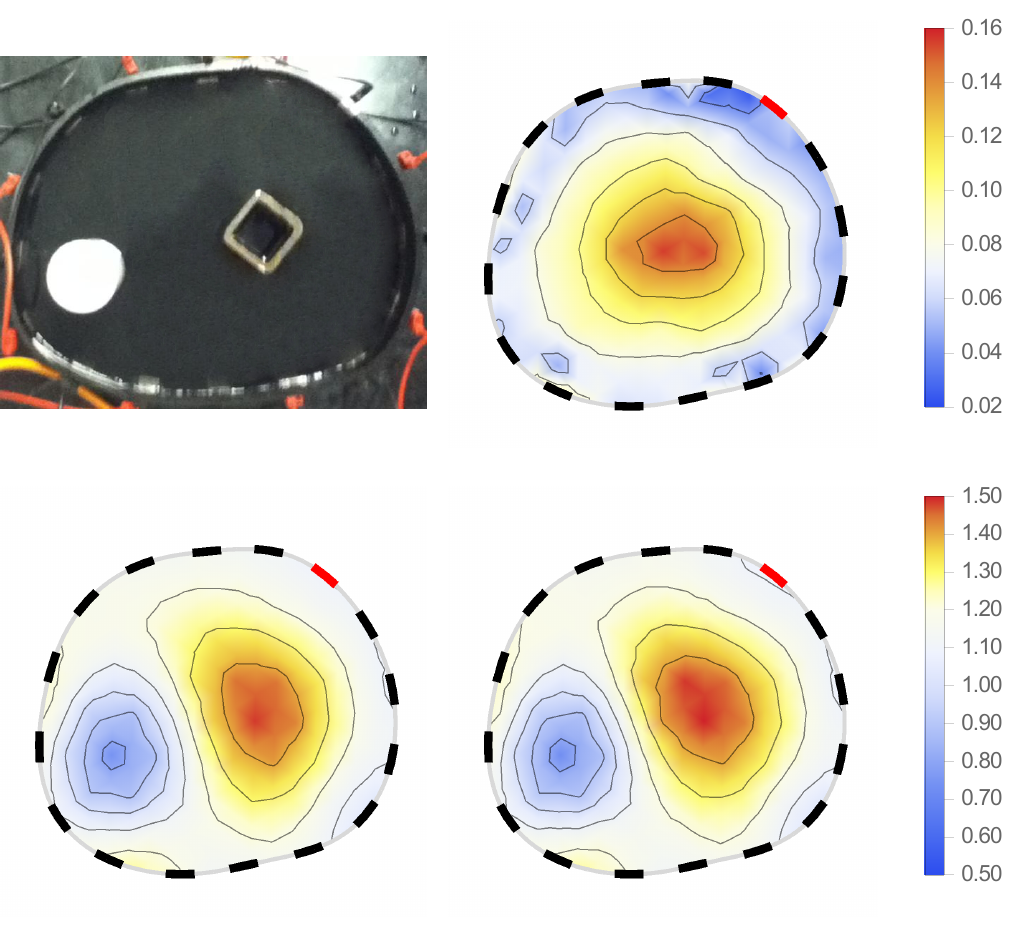}}
  \end{center}
  \caption{
    Results of the fifth example.
    Top left: the target with one insulating and one highly conducting inclusion.
    Top right: the SD estimate $\sigma_{\rm SD}$.
    Bottom left: the MAP estimate $\sigma_{\rm MAP}$.
    Bottom right: the CM estimate $\sigma_{\rm CM}$.
    The unit in all images is mS, the MAP and CM estimates use the same colormap, and the colorbar tick markers correspond to the contour lines in the images.
}
  \label{fig:Example1_5}
\end{figure}

We conclude with two experiments with a pair of embedded inclusions: one plastic and one metallic cylinder. The target configurations are shown in the top left images of Figures~\ref{fig:Example1_4} and \ref{fig:Example1_5}. The other images in Figures~\ref{fig:Example1_4} and \ref{fig:Example1_5} illustrate the corresponding MAP, CM, and SD estimates for the two measurement configurations. Even in this slightly more complicated setting, our algorithm produces reasonably good reconstructions: in both experiments, the positions of the two inhomogeneities can be identified accurately from the MAP and CM estimates. Indeed, the highest and lowest
reconstructed conductivity levels are attained close to the center points of the metallic and plastic inclusions, respectively. However, the reconstructions are heavily blurred, which is not very surprising as the employed prior \eqref{prior2} prefers slow changes over sharp boundaries. 

We have not tested the algorithm with a higher number of inhomogeneities, but we suspect that the parametrization of the conductivity by the $76$ pixels depicted in Figure~\ref{fig:pixelgrid} is insufficient for reconstructing much more complicated phantoms than the ones in Figures~\ref{fig:Example1_4} and \ref{fig:Example1_5}.

\section{Conclusions}
\label{sec:conclusions}

We have studied the feasibility of solving the reconstruction problem of EIT by combining SCEM and sGFEM, with the unknown conductivity field parametrized by its values at a set of pixels.
The functionality of the method was demonstrated by applying it to five data sets from water tank experiments. In all cases, the resulting MAP and CM estimates clearly provided useful information about the conductivity phantom.

Assuming that the measurement configuration and the preliminary bounds for the pixelwise conductivity values are known well in advance, the pre-measurement phase of the reconstruction algorithm can be performed off-line, and subsequently
the (approximate) posterior distribution of the conductivity is obtained practically for free when the measurement data becomes available. Hence, the on-line solution phase of the algorithm consists solely of extracting the desired estimators from the explicitly parametrized posterior.

In the post-measurement phase of the algorithm, we resorted exclusively to a Gaussian prior with a covariance matrix of the type \eqref{equ:priorcov}, which resulted in blurred conductivity reconstructions. In principle, it should also be possible to use any other prior (e.g., total variation \cite{Vogel96}) for the conductivity in the post-processing phase. Such a modification would only affect the form of the target function in \eqref{equ:MAPproblem} and the integrands in \eqref{CMint1} and \eqref{CMint5}, but it could lead to, e.g., more accurate detection of inclusion boundaries. This line of research is left for future studies.

\section*{Acknowledgments}
We would like to thank Professor Jari Kaipio's research group at the University of Eastern Finland (Kuopio) for granting us access to their EIT devices. 
We acknowledge CSC -- IT Center for Science Ltd. for the allocation of computational resources (project ay6302).

\bibliographystyle{acm}
\bibliography{sGFEM_with_local_conductivity_basis_for_EIT}

\end{document}